\documentclass{gtmon_a}
\pdfoutput=1

\proceedingstitle{The Zieschang Gedenkschrift}
\conferencestart{5 September 2007}
\conferenceend{8 September 2007}
\conferencename{Conference in honour of Heiner Zieschang}
\conferencelocation{Toulouse, France}

\editor{Michel Boileau}
\givenname{Michel}
\surname{Boileau}

\editor{Martin Scharlemann}
\givenname{Martin}
\surname{Scharlemann}

\editor{Richard Weidmann}
\givenname{Richard}
\surname{Weidmann}

\title{On the number of optimal surfaces}

\author{Alina Vdovina}
\givenname{Alina}
\surname{Vdovina}
\address{School of Mathematics and Statistics\\Newcastle University\\\newline
Newcastle-upon-Tyne NE1 7RU\\UK}
\email{alina.vdovina@ncl.ac.uk}
\urladdr{}

\dedicatory{To the memory of Heiner Zieschang}

\volumenumber{14}
\issuenumber{}
\publicationyear{2008}
\papernumber{24}
\startpage{557}
\endpage{567}

\doi{}
\MR{}
\Zbl{}


\received{31 July 2006}
\revised{28 July 2007}
\accepted{28 July 2007}
\published{29 April 2008}
\publishedonline{29 April 2008}
\proposed{}
\seconded{}
\corresponding{}
\version{}

\arxivreference{}

\makeatletter
\def\cnewtheorem#1[#2]#3{\newtheorem{#1}{#3}[section]
\expandafter\let\csname c@#1\endcsname\c@Theorem}
\makeatother

\AtBeginDocument{\let\bar\wbar\let\tilde\wtilde\let\hat\what}
\def\strutt{\vrule width 0pt height 12pt}
\newcommand{\N}{\mathbb N}

\newtheorem{Theorem}{Theorem}[section]
\cnewtheorem{Corollary}[Theorem]{Corollary}
\cnewtheorem{Proposition}[Theorem]{Proposition}
\theoremstyle{definition}
\cnewtheorem{Definition}[Theorem]{Definition}

\makeautorefname{Theorem}{Theorem}

\begin{document}

\begin{asciiabstract}
Let X be a closed oriented Riemann surface of genus > 1 of constant
negative curvature -1.  A surface containing a disk of maximal radius
is an optimal surface.  This paper gives exact formulae for the number
of optimal surfaces of genus > 3 up to orientation-preserving
isometry.  We show that the automorphism group of such a surface is
always cyclic of order 1,2,3 or 6.  We also describe a combinatorial
structure of nonorientable hyperbolic optimal surfaces.

\end{asciiabstract}

\begin{htmlabstract}
Let X be a closed oriented Riemann surface of genus &ge; 2 of
constant negative curvature -1.  A surface containing a disk of
maximal radius is an <em>optimal surface</em>.  This paper gives exact
formulae for the number of optimal surfaces of genus &ge; 4 up to
orientation-preserving isometry.  We show that the automorphism group
of such a surface is always cyclic of order 1, 2, 3 or 6. We also
describe a combinatorial structure of nonorientable hyperbolic optimal
surfaces.
\end{htmlabstract}

\begin{abstract}
Let $X$ be a closed oriented Riemann surface of genus $\geq 2$ of
constant negative curvature $-1$.  A surface containing a disk of
maximal radius is an {\em optimal surface}.  This paper gives exact
formulae for the number of optimal surfaces of genus $\geq 4$ up to
orientation-preserving isometry.  We show that the automorphism group
of such a surface is always cyclic of order 1, 2, 3 or 6. We also
describe a combinatorial structure of nonorientable hyperbolic optimal
surfaces.
\end{abstract}

\maketitle

\section*{Introduction}

 Let $X$ be a compact Riemann surface of genus $\geq 2$
of constant negative curvature $-1$.
We consider the maximal radius of an embedded open metric
disk in $X$.
 A surface containing such 
disk is a {\em optimal surface}.

Such surfaces, obtained from generic polygon side pairings, appear in
the literature in different contexts and go back to Fricke and Klein
(see Girondo and Gonz\'alez-Diez \cite{[GG2]} for a great survey).

The radius $R_g$ of an maximal
 embedded disk, as well as the radius 
 $C_g$ of an minimal covering disk were computed by Bavard in \cite{[Ba]}:
 $$R_g=\cosh^{-1} \frac{1}{2\sin \beta_g}, \beta_g=\pi/{(12g-6)}.$$
  $$C_g=\cosh^{-1}\frac{1}{\sqrt3 \tan \beta_g}, \beta_g=\pi/{(12g-6)}.$$
 The discs of maximal radius occur
in those surfaces which admit as Dirichlet domain 
a regular polygon with the largest possible number of
sides $12g-6$.

We give an exact formula for the number of optimal surfaces of genus
$g\geq 4$, up to orientation-preserving isometry, as well as an
explicit construction of all optimal surfaces of genus $g\geq 4$. Note
that, in this paper, we consider surfaces always up to
orientation-preserving isometries.  We show, that for genus $g \geq 4$
the automorphism group of an oriented optimal surface is always cyclic
of order 1, 2, 3 or 6 and we give an explicit formula for the number
of nonisometric optimal surfaces.  It follows from the formula, that
the number grows factorially with $g$ (more precisely, it grows as
$(2g)!$). This is a significant improvement compared to \cite{[GG2]},
where it was noted that the number grows exponentially with $g$.

Also we give explicit formulae of optimal surfaces having exactly $d$
automorphisms, where $d$ is 1, 2, 3 or 6.  These formulae show that
asymptotically almost all optimal surfaces have no automorphisms.  In
particular, for $d=1$ we have a big family of explicitly constructed
surfaces with no automorphisms.  Let us note that another families of
surfaces with no automorphisms were considered by Everitt in
\cite{[E]}, Turbek in \cite{[T]} and by Girondo and Gonz\'alez-Diez
in \cite{[GG2]}.  The questions of explicit construction, enumeration
and description of automorphisms of genus $2$ optimal surfaces were
solved by Girondo and Gonz\'alez-Diez in \cite{[GG3]}.  C Bavard
\cite{[Ba]} proved, that if a surface contains an embedded disk of
maximal radius if and only if it admits a covering disk of minimal
radius and that optimal surfaces are modular curves.

 We will show, that oriented maximal
 Wicks forms
 and optimal surfaces
are in bijection for $g \geq 4$.

Wicks forms are canonical forms for products of commutators in free
groups (Vdovina \cite{[V]}). Wicks forms arise as well in a much
broader context of connection of branch coverings of compact surfaces
and quadratic equations in a free group (Bogatyi, Gon{\c{c}}alves,
Kudryavtseva, Weidmann and Zieschang \cite{6,5,7}).
But in the present paper we restrict ourselves to a particular type of
Wicks forms related to optimal surfaces.

For $g=2$ the bijection between Wicks forms and optimal surfaces was
proved by Girondo and Gonz\'alez-Diez in \cite{[GG1]} and for $g=3$ the question is still open.

\fullref{sec:1} formulates our main results and 
introduces oriented Wicks forms (cellular
 decompositions with only one face of oriented surfaces).
We include the detailed explanation of results presented
in Bacher and Vdovina \cite{[BV]} for completeness.

\fullref{sec:2} contains the proof of our main results. In 
\fullref{sec:3} we treat the case of nonorientable surfaces.

\section{Main results}
\label{sec:1}

\begin{Definition}\label{def:1.1}
 An {\it oriented Wicks form\/} is a cyclic word $w= w_1w_2\dots w_{2l}$
 (a cyclic word is the orbit of a linear word under cyclic permutations)
 in some alphabet\break $a_1^{\pm 1},a_2^{\pm 1},\dots, a_l^{\pm}$ of letters
 $a_1,a_2,\dots, a_l$ and their inverses $a_1^{-1},a_2^{-1},\dots, a_l^{-1}$
 such that
\begin{itemize}
\item[(i)] if $a_i^\epsilon$ appears in $w$ (for $\epsilon\in\{\pm 1\}$)
 then $a_i^{-\epsilon}$ appears exactly once in $w$,
\item[(ii)] the word $w$ contains no cyclic factor (subword of
 cyclically consecutive letters in $w$) of the form $a_i a_i^{-1}$ or
 $a_i^{-1}a_i$ (no cancellation),
\item[(iii)] if $a_i^\epsilon a_j^\delta$ is a cyclic factor of $w$ then
 $a_j^{-\delta}a_i^{-\epsilon}$ is not a cyclic factor of $w$ 
(substitutions of the form
 $a_i^\epsilon a_j^\delta\longmapsto x,
 \quad a_j^{-\delta}a_i^{-\epsilon}\longmapsto x^{-1}$ are impossible).
\end{itemize}
\end{Definition}

 An oriented Wicks form $w=w_1w_2\dots$ in an alphabet $A$ is
 {\em isomorphic\/} to $w'=w'_1w'_2$ in an alphabet $A'$ if
 there exists a bijection $\varphi:A\longrightarrow A'$ with
 $\varphi(a^{-1})=\varphi(a)^{-1}$ such that $w'$ and
 $\varphi(w)=\varphi(w_1)\varphi(w_2)\dots$ define the
 same cyclic word.

 An oriented Wicks form $w$ is an element of the commutator subgroup
 when considered as an element in the free group $G$ generated by
 $a_1,a_2,\dots$. We define the {\em algebraic genus\/} $g_a(w)$ of
 $w$ as the least positive integer $g_a$ such that $w$ is a product
 of $g_a$ commutators in $G$.

 The {\em topological genus\/} $g_t(w)$ of an oriented Wicks
 form $w=w_1\dots w_{2e-1}w_{2e}$ is defined as the topological
 genus of the oriented compact connected surface obtained by
 labelling and orienting the edges of a $2e$--gon (which we
 consider as a subset of the oriented plane) according to
 $w$ and by identifying the edges in the obvious way.

\begin{Proposition}[Culler \cite{[C]} and Comerford and Edmunds
\cite{[CE]}]\label{prop:1.1}
The algebraic genus and the topological genus of an oriented Wicks
 form coincide.\end{Proposition}

 We define the {\em genus\/} $g(w)$ of an oriented
 Wicks form $w$ by $g(w)=g_a(w)=g_t(w)$.

 Consider the oriented compact surface $S$ associated to an oriented 
 Wicks form $w=w_1\dots w_{2e}$. This surface carries an embedded graph
 $\Gamma\subset S$ such that $S\setminus \Gamma$ is an open polygon
 with $2e$ sides (and hence connected and simply connected).
 Moreover, conditions (ii) and (iii) on Wicks form imply that $\Gamma$ 
 contains no vertices of degree $1$ or $2$ (or equivalently that the
 dual graph of $\Gamma\subset S$ contains no faces which are $1$--gons
 or $2$--gons). This construction works also
 in the opposite direction: Given a graph $\Gamma\subset S$
 with $e$ edges on an oriented compact connected surface $S$ of genus $g$
 such that $S\setminus \Gamma$ is connected and simply connected, we get
 an oriented Wicks form of genus $g$ and length $2e$ by labelling and 
 orienting the edges of $\Gamma$ and by cutting $S$ open along the graph
 $\Gamma$. The associated oriented Wicks form is defined as the word
 which appears in this way on the boundary of the resulting polygon
 with $2e$ sides. We identify henceforth oriented Wicks
 forms with the associated embedded graphs $\Gamma\subset S$, when
 speaking of vertices and edges of an oriented Wicks form.

 The formula for the Euler characteristic
 $$\chi(S)=2-2g=v-e+1$$
 (where $v$ denotes the number of vertices and $e$ the number
 of edges in $\Gamma\subset S$) shows that
 an oriented Wicks form of genus $g$ has at least length $4g$
 (the associated graph has then a unique vertex of degree $4g$
 and $2g$ edges) and at most length $6(2g-1)$ (the associated
 graph has then $2(2g-1)$ vertices of degree three and
 $3(2g-1)$ edges).

 We call an oriented Wicks form of genus $g$ {\em maximal\/} if it has
 length $6(2g-1)$. Oriented maximal Wicks forms are dual to 1--vertex
 triangulations. This can be seen by cutting the oriented surface $S$ 
 along $\Gamma$, hence obtaining a polygon $P$ with $2e$ sides. 
 We draw a star $T$ on $P$ which joins an interior point of $P$
 with the midpoints of all its sides. Regluing $P$ we recover $S$
 which carries now a 1--vertex triangulation given by $T$ and each
 1--vertex triangulation is of this form for some oriented maximal
 Wicks form
 
 This construction shows that we can work indifferently with
 1--vertex triangulations or with oriented maximal Wicks forms.

 A vertex $V$ of degree three 
(with oriented edges $a,b,c$ pointing toward $V$) is
 {\em positive\/} if
$$w=ab^{-1}\dots bc^{-1}\dots ca^{-1}\dots \quad {\rm or }\quad
 w=ac^{-1}\dots cb^{-1}\dots ba^{-1}\dots $$
 and $V$ is {\em negative\/} if    
 $$w=ab^{-1}\dots ca^{-1}\dots bc^{-1}\dots \quad {\rm or }\quad 
 w=ac^{-1}\dots ba^{-1}\dots ab^{-1}\dots \ .$$
 The {\em automorphism group\/} ${\rm Aut}(w)$ of an oriented
 Wicks form $$w=w_1w_2\dots w_{2e}$$ of length $2e$ is the group of all
 cyclic permutations $\mu$ of the linear word $w_1w_2\dots w_{2e}$ such
 that $w$ and $\mu(w)$ are isomorphic linear words (ie, $\mu(w)$ is
 obtained from $w$ by permuting the letters of the alphabet). The group
 ${\rm Aut}(w)$ is a subgroup of the cyclic group ${\Z}/2e{\Z}$
 acting by cyclic permutations on linear words representing $w$.

 The automorphism group ${\rm Aut}(w)$ of an oriented Wicks
 form can of course also be described in terms of permutations on the
 oriented edge set induced by orientation-preserving
 homeomorphisms of $S$ leaving $\Gamma$ invariant. In particular
 an oriented maximal Wicks form and the associated dual
 1--vertex triangulation have isomorphic automorphism groups.

 We define the {\em mass\/} $m(W)$ of a finite set $W$ of oriented
 Wicks forms by
 $$m(W)=\sum_{w\in W}{1\over \vert{\rm Aut}(w)\vert}\ .$$
 Let us introduce the sets:
\par $W_1^g$: all oriented maximal
 Wicks forms of genus $g$ (up to isomorphism);
\par $W^g_2(r)\subset W_1^g$: all oriented 
 maximal Wicks forms having an automorphism of order $2$ leaving 
 exactly $r$ edges of $w$ invariant by reversing their
 orientation (this automorphism is the half-turn with respect to the
 \lq\lq midpoints" of these edges and exchanges the two adjacent
 vertices of an invariant edge);
\par $W^g_3(s,t)\subset W_1^g$: all oriented maximal 
 Wicks forms having an automorphism of order $3$ leaving
 exactly $s$ positive and $t$ negative vertices invariant
 (this automorphism permutes cyclically the edges around
 an invariant vertex);
\par $W^g_6(3r;2s,2t)=W^g_2(3r)\cap W^g_3(2s,2t)$: all oriented
 maximal Wicks forms having an automorphism $\gamma$ of order 6
 with $\gamma^3$ leaving $3r$ edges invariant and $\gamma^2$
 leaving $2s$ positive and $2t$ negative vertices invariant
 (it is useless to consider the set $W_6^g(r';s',t')$ defined
 analogously since $3$ divides $r'$ and $2$ divides $s',t'$ if
 $W_6^g(r';s',t')\not=\emptyset$).

\par We define now the {\em masses\/} of these sets as
$$\begin{matrix}
   m_1^g\hfill&=&\displaystyle \sum_{w\in W_1^g}
    {1\over \vert{\rm Aut}(w)\vert}\ ,\hfill \\
    m_2^g(r)\hfill&=&\displaystyle \sum_{w\in W_2^g(r)}
    {1\over \vert{\rm Aut}(w)\vert}\ ,\hfill \\
    m_3^g(s,t)\hfill&=&\displaystyle \sum_{w\in W_3^g(s,t)}
    {1\over \vert{\rm Aut}(w)\vert}\ ,\hfill \\
    m_6^g(3r;2s,2t)\hfill&=&\displaystyle\sum_{w\in W_6^g(3r;2s,2t)}
    {1\over \vert{\rm Aut}(w)\vert}\ .\hfill 
\end{matrix}
$$

\begin{Theorem}{\rm\cite{[BV]}}\label{Theo:1.1}\

\begin{itemize}
\item[\rm(i)] The group ${\rm Aut}(w)$ of automorphisms of an oriented
 maximal Wicks form $w$ is cyclic of order $1,\ 2,\ 3$ or $6$.

\item[\rm(ii)] $\displaystyle
 m_1^g={2\over 1}\Big({1^2\over 12}\Big)^g{(6g-5)!\over g!(3g-3)!}
 \ .$

\item[\rm(iii)] $m_2^g(r)>0$ (with $r\in {\N}$) if and only if
 $f={2g+1-r\over 4}\in \{0,1,2,\dots\}$ and we have then
\par $\displaystyle m^g_2(r)={2\over 2}\Big({2^2\over 12}\Big)^f
 {1\over r!}{(6f+2r-5)!\over f!(3f+r-3)!}
 \ .$

\item[\rm(iv)] $m_3^g(s,t)>0$(with $r,s \in {\N}$) if and only if
 $f={g+1-s-t\over 3}\in \{0,1,2,\dots\}$, $s\equiv 2g+1\pmod 3$ and
 $t\equiv 2g\pmod 3$ (which follows from the two previous conditions).
 We have then
\par $\displaystyle m^g_3(s,t)={2\over 3}\Big({3^2\over 12}\Big)^f
 {1\over s!t!}{(6f+2s+2t-5)!\over f!(3f+s+t-3)!}$ if $g>1$ and 
$\displaystyle m^1_3(0,2)={1\over 6}$.

\item[\rm(v)] $m_6^g(3r;2s,2t)>0$(with $r,s,t \in {\N}$) if and only if
 $f={2g+5-3r-4s-4t \over 12}\in \{0,1,2,\dots\}$,
 $2s\equiv 2g+1\pmod 3$ and $2t\equiv 2g\pmod 3$ (follows in fact
 from the previous conditions). We have then
\par $\displaystyle
 m_6^g(3r;2s,2t)={2\over 6}\Big({6^2\over 12}\Big)^f
 {1\over r!s!t!}{(6f+2r+2s+2t-5)!\over f!(3f+r+s+t-3)!}$\\
if $g>1$ and $\displaystyle m^1_6(3;0,2)={1\over 6}$.
\end{itemize}
\end{Theorem}

\begin{Theorem}\label{Theo:1.2}\

\begin{itemize}
\item[\rm(i)] The group ${\rm Aut}(S_g)$ of automorphisms of an optimal
 surface is cyclic of order $1,\ 2,\ 3$ or $6$ in every genus $\geq 4$.

\item[\rm(ii)] There is a bijection between  isomorphism
 classes of oriented maximal genus $g$ Wicks forms
and optimal genus $g$ surfaces  for $g \geq 4$ up to orientation-preserving
isometries.
\end{itemize}
\end{Theorem}

Set
 $$\begin{matrix}
m_2^g&=&\sum_{r\in {\N},\ (2g+1-r)/4\in {\N}\cup \{0\}}
 m_2^g(r)\ ,\hfill \\
\strutt m_3^g&=&\sum_{s,t\in {\N},\ (g+1-s-t)/3\in {\N}\cup \{0\},\ 
 s\equiv 2g+1 ({\rm mod} 3)} m_3^g(s,t)\ ,\hfill \\
\strutt m_6^g&=&\sum_{r,s,t\in {\N},\ (2g+5-3r-4s-4t)/12\in
 {\N}\cup \{0\},\ 2s \equiv 2g+1({\rm mod} 3)} m_6^g(3r;2s,2t)\hfill
\end{matrix}$$
 (all sums are finite) and denote by $M_d^g$ the number of automorphism
 classes of optimal genus $g$ surfaces having an
  automorphism of order $d$ (ie, an automorphism group with order
 divisible by $d$, see \fullref{Theo:1.2}(i)).

Combining \fullref{Theo:1.2}(ii) and \cite[Theorem 1.3]{[BV]},
we obtain the following:

\begin{Corollary}\label{Cor:1.1}
For $\geq 4$, we have
 $$\begin{matrix}
M_1^g=m_1^g+m_2^g+2m_3^g+2m_6^g\ ,\hfill \\
\strutt M_2^g=2m_2^g+4m_6^g\ ,\hfill \\
\strutt M_3^g=3m_3^g+3m_6^g\ ,\hfill \\
\strutt M_6^g=6m_6^g \hfill
\end{matrix} $$
 and $M_d^g=0$ if $d$ is not a divisor of $6$.
\end{Corollary}

The number $M_1^g$ in this Theorem is  the number of
 optimal surfaces
 of genus $g$ for $g \geq 4$ up to orientation-preserving
isometry. The first 14 values $M_1^2,\dots,M_1^{15}$, except $M_1^{3}$, are
displayed in \fullref{Table1}.

The following result is an immediate consequence of \fullref{Theo:1.2}(i).

\begin{Corollary}\label{Cor:1.2}For $g \geq 4$ there are exactly
\par $M_6^g$ nonisometric optimal surfaces with $6$
 automorphisms,\newline\strutt $M_3^g-M_6^g$ nonisometric optimal
 surfaces with $3$ automorphisms,\newline\strutt $M_2^g-M_6^g$
 nonisometric optimal surfaces with $2$ automorphisms
 and\newline\strutt $M_1^g-M_2^g-M_3^g+M_6^g$ nonisometric optimal
 surfaces without non-trivial automorphisms.
\end{Corollary}

\begin{table}[ht!]\small
\caption{The number of oriented optimal surfaces of genus $2,4 \dots,15$}
\label{Table1}
\centerline{$\begin{matrix}
%1\hfill &  1\hfill \\ 
2\hfill &  9\hfill \\ 
%3\hfill &  1726\hfill \\ 
4\hfill &  1349005\hfill \\ 
5\hfill &  2169056374\hfill \\ 
6\hfill &  5849686966988\hfill \\ 
7\hfill &  23808202021448662\hfill \\ 
8\hfill &  136415042681045401661\hfill \\ 
9\hfill &  1047212810636411989605202\hfill \\ 
10\hfill &  10378926166167927379808819918\hfill \\ 
11\hfill &  129040245485216017874985276329588\hfill \\ 
12\hfill &  1966895941808403901421322270340417352\hfill \\ 
13\hfill &  36072568973390464496963227953956789552404\hfill \\ 
14\hfill &  783676560946907841153290887110277871996495020\hfill \\ 
15\hfill &  19903817294929565349602352185144632327980494486370\hfill \\ 
\end{matrix}$}
\end{table}

\section[Proof of \ref{Theo:1.2}]{Proof of \fullref{Theo:1.2}}
\label{sec:2}

{\bf Proof of (i)}

Let $w$ be an oriented maximal Wicks form with an automorphism
$\mu$ of order $d$. Let $p$ be a prime dividing $d$. The automorphism
$\mu'=\mu^{d/p}$ is hence of order $p$. If $p\not= 3$ then $\mu'$
acts without fixed vertices on $w$ and \cite[Proposition 2.1]{[BV]}
 shows that $p$
divides the integers $2(g-1)$ and $2g$ which implies $p=2$. 
The order $d$ of $\mu$ is hence 
of the form $d=2^a3^b$. Repeating the above argument with the prime 
power $p=4$ shows
that $a\leq 1$. 

\par All orbits of $\mu^{2^a}$ on the set
of positive (respectively negative)
vertices have either $3^b$ or $3^{b-1}$ elements and this leads 
to a contradiction if $b\geq 2$. 
This shows that $d$ divides $6$ and proves that 
the automorphism groups of oriented maximal
 Wicks forms
are always cyclic of order 1, 2, 3 or 6.

Let us consider an optimal genus $g$ surface $S_g$.
It was proved in \cite{[Ba]}, that a surface is optimal
if and only if it can be obtained from a regular oriented hyperbolic
$(12g-6)$--gon  with angles $2\pi/3$ such that the image
of the boundary of the polygon after identification
of corresponded sides is a geodesic graph with $4g-2$ vertices
of degree $3$ and $6g-3$ edges of equal length.
It was shown in $\cite{[GG3]}, g \geq 4$, that any isometry of
an optimal surface of genus $g>3$ is realized by
a rotation of the $(12g-6)$--gon.  

Let $P$ be a regular geodesic hyperbolic polygon with $12g-6$ equal sides
and all angles equal to $2\pi/3$, equipped with a oriented maximal
genus $g$ Wicks form $W$ on its boundary. Consider the surface $S_g$
obtained from $P$ by identification of sides with the same labels. 
Since we made the identification using an oriented maximal
 Wicks form of length
$12g-6$, the boundary of $D$ becomes a graph $G$
 with $4g-2$ vertices
of degree $3$ and $6g-3$ edges (see \fullref{sec:1}).  
We started from a regular geodesic hyperbolic polygon with angles
$2\pi/3$, so $G$ is a geodesic graph with edges of equal length
and the surface
$S_g$ is optimal.

So, the surface is optimal if and only if it
can be obtained from a regular hyperbolic polygon with $12g-6$ equal sides
and all angles equal to $2\pi/3$, equipped with an oriented maximal
genus $g$ Wicks form $W$ on its boundary.
The isometry of $S_g$ must be realized by a rotation
of the $(12g-6)$--gon \cite{[GG3]}, so the isometry must be an automorphism
of the Wicks form. Since the automorphism groups of oriented maximal
 Wicks forms
are always cyclic of order 1, 2, 3 or 6, the automorphism
groups of genus $g \geq 4$ optimal surfaces are also
 cyclic of order 1, 2, 3 or 6.

\medskip
{\bf Proof of (ii)}

Every oriented maximal Wicks form defines exactly one oriented optimal
surface, namely the surface obtained from 
 a regular hyperbolic polygon with $12g-6$ equal sides
and all angles equal to $2\pi/3$ with an oriented maximal
genus $g$ oriented maximal Wicks form $W$ on its boundary.

So, to prove the bijection between the set of equivalence
 classes of oriented maximal genus $g$ Wicks forms
and optimal genus $g$ surfaces for $g \geq 4$ we need
to show, that for every optimal surface $S_g$ there is only
one oriented maximal Wicks form $W$ such that $S_g$
can be obtained from a regular hyperbolic $12g-6$--gon
with $W$ on its boundary. It was shown in \cite{[GG3]}, that
for $g \geq 4$ the  maximal open disk $D$ of radius 
 $R_g=\cosh^{-1}(1/2\sin \beta_g)$, $\beta_g=\pi/{(12g-6)}$,
 embedded in $S_g$ is unique. 
Consider the center $c$ of the disk $D$. The discs of radius
$R_g$ with the centers in the images of $c$
in the universal covering of $S_g$ form a packing of
the hyperbolic plane by discs.
To this packing one can classically associate a tessellation $T$
of the hyperbolic plane by regular $12-6$--gons, which
are Dirichlet domains for $S_g$. And such a tessellation
is unique because of negative curvature.
The fundamental group of the surface $S_g$ naturally acts on
the hyperbolic plane by covering transformations preserving the
tesselation $T$. But each such action defines a Wicks form.
\fullref{Theo:1.2} is proved. \qed

\section{Nonorientable optimal surfaces}
\label{sec:3}

In a similar to the orientable case way
we can associate nonorientable optimal surfaces with
nonoriented Wicks forms.

\begin{Definition}\label{Def:3.1}
 A {\it nonoriented Wicks form\/} is a cyclic word $w= w_1w_2\dots w_{2l}$
 in some alphabet $a_1^{\pm 1},a_2^{\pm 1},\dots, a_l^{\pm}$ of letters
 $a_1,a_2,\dots,a_l$ and their inverses $a_1^{-1},a_2^{-1},$ $\dots,a_l^{-1}$
 such that
\begin{itemize}
\item[(i)] every letter appears exactly twice in $w$, and at least
one letter appears with the same exponent,
\item[(ii)] the word $w$ contains no cyclic factor (subword of
 cyclically consecutive letters in $w$) of the form $a_i a_i^{-1}$ or
 $a_i^{-1}a_i$ (no cancellation),
\item[(iii)] if $a_i^\epsilon a_j^\delta$ is a cyclic factor of $w$ then
 $a_j^{-\delta}a_i^{-\epsilon}$ is not a cyclic factor of $w$ 
(substitutions of the form
 $a_i^\epsilon a_j^\delta\longmapsto x,
 \quad a_j^{-\delta}a_i^{-\epsilon}\longmapsto x^{-1}$ are impossible).
\end{itemize}
\end{Definition}

 A nonoriented Wicks form $w=w_1w_2\dots$ in an alphabet $A$ is
 {\em isomorphic\/} to $w'=w'_1w'_2$ in an alphabet $A'$ if
 there exists a bijection $\varphi:A\longrightarrow A'$ with
 $\varphi(a^{-1})=\varphi(a)^{-1}$ such that $w'$ and
 $\varphi(w)=\varphi(w_1)\varphi(w_2)\dots$ define the
 same cyclic word.

 A nonoriented Wicks form $w$ is an element of the subgroup,
consisting of products of squares,
 when considered as an element in the free group $G$ generated by
 $a_1,a_2,\dots$. We define the {\em algebraic genus\/} $g_a(w)$ of
 $w$ as the least positive integer $g_a$ such that $w$ is a product
 of $g_a$ squares in $G$.

 The {\em topological genus\/} $g_t(w)$ of a nonoriented Wicks
 form $w=w_1\dots w_{2e-1}w_{2e}$ is defined as the topological
 genus of the nonorientable compact connected surface obtained by
 labelling and orienting the edges of a $2e$--gon  according to
 $w$ and by identifying the edges in the obvious way.
The algebraic and the topological genus of a nonoriented Wicks
 form coincide (cf \cite{[CE],[C]}).

 We define the {\em genus\/} $g(w)$ of a nonoriented
 Wicks form $w$ by $g(w)=g_a(w)=g_t(w)$.

 Consider the nonorientable compact surface $S$ associated to a nonoriented 
 Wicks form $w=w_1\dots w_{2e}$. This surface carries an immersed graph
 $\Gamma\subset S$ such that $S\setminus \Gamma$ is an open polygon
 with $2e$ sides (and hence connected and simply connected).
 Moreover, conditions (ii) and (iii) on Wicks form imply that $\Gamma$ 
 contains no vertices of degree $1$ or $2$ (or equivalently that the
 dual graph of $\Gamma\subset S$ contains no faces which are $1$--gons
 or $2$--gons). This construction works also
 in the opposite direction: Given a graph $\Gamma\subset S$
 with $e$ edges on a nonorientable compact connected surface $S$ of genus $g$
 such that $S\setminus \Gamma$ is connected and simply connected, we get
 a nonoriented Wicks form of genus $g$ and length $2e$ by labelling and 
 orienting the edges of $\Gamma$ and by cutting $S$ open along the graph
 $\Gamma$. The associated nonoriented Wicks form is defined as the word
 which appears in this way on the boundary of the resulting polygon
 with $2e$ sides. We identify henceforth nonoriented Wicks
 forms with the associated immersed graphs $\Gamma\subset S$,
 speaking of vertices and edges of nonorientable Wicks form.

 The formula for the Euler characteristic
 $$\chi(S)=2-g=v-e+1$$
 (where $v$ denotes the number of vertices and $e$ the number
 of edges in $\Gamma\subset S$) shows that
 a nonoriented Wicks form of genus $g$ has at least length $2g$
 (the associated graph has then a unique vertex of degree $2g$
 and $g$ edges) and at most length $6(g-1)$ (the associated
 graph has then $2(g-1)$ vertices of degree three and
 $3(g-1)$ edges).

 We call a nonoriented Wicks form of genus $g$ {\em maximal\/} if it has
 length $6(g-1)$.

It follows from \cite{[Ba]}, that any nonorientable optimal
surface of genus $g\geq 3$ can be obtained from
a regular hyperbolic $6g-6$--gon, with angles $2\pi/3$,
with a nonoriented Wicks form of genus $g$ on its boundary,
by identification of corresponded sides respecting
orientation. So, we have that the number of nonorientable
genus $g$ optimal surfaces is majorated by the number $M(g)$ of
nonorientable genus $g$ Wicks forms. The asymptotic behaviour
of $M(g)$ is established in \cite{[A1]}. The description
and classification of nonorientable Wicks forms is done in \cite{[A2]}.

\bibliographystyle{gtart}
\bibliography{link}

\end{document}